\documentclass[11pt]{article}
\usepackage{amsmath,amssymb,amsthm}

\newcommand{\res}{\operatorname{Res}}
\newcommand{\coker}{\operatorname{coker}}
\newcommand{\caB}{\mathcal B}
\newcommand{\caC}{\mathcal C}
\newcommand{\caD}{\mathcal D}
\newcommand{\omCB}{\omega_{\caC/\caB}}

\newtheorem{lemma}{Lemma}
\newtheorem{defn}[lemma]{Definition}

\usepackage{graphicx}
\graphicspath{%
    {converted_graphics/}
    {CBMS/CBMSbook/}
}
\begin{document}

\title{On families of holomorphic differentials on degenerating annuli}         
\author{Scott A. Wolpert\footnote{Partially supported by National Science Foundation grant DMS - 1005852.}}        
\date{November 23, 2011}          
\maketitle
\vspace{-.3in}
\begin{center}
{\em To Cliff Earle, leading a way by example.}
\end{center}

\section{Introduction.}
We consider the local analytic behavior for a family of holomorphic differentials on a family of degenerating annuli.  The matter is closely related to the situation for a neighborhood of a node in a standard nodal family of Riemann surfaces/curves.  Basic considerations for a family of differentials include the local analytic description on annuli cores, the role of the relative dualizing sheaf as a setting for a family, extension properties of a family on the limiting nodal space and divisors of families. 

We present three results and discussion.  The first is the normal families formulation Lemma \ref{normal}.  The second is an isomorphism of sheaves, formula (\ref{psiPsi}), giving a direct description of families of regular $k$-differentials (sections of powers of the relative dualizing sheaf) in terms of $k$-canonical forms on the total space of the family. The third is a general holomorphic extension property, Lemma \ref{extension}, for families given on smooth Riemann surfaces/curves to extend to the limiting nodal Riemann surfaces/curves. 

The module of an annulus/conformal cylinder is defined as the reciprocal extremal length for the family of rectifiable curves separating boundaries \cite[Chap. 1]{Ahqc}.  For a geometric annulus $A=\{r_1<|\zeta|<r_2\}$ in $\mathbb C$, the module is given as $\mathbf M(A)=(\log r_2/r_1)/2\pi$.  Basic for a geometric annulus are the rotationally invariant vector field $\zeta\partial/\partial \zeta$ and dual differential $d\zeta/\zeta$.  The geometric annulus $A$ is equivalent to the Euclidean cylinder $\{0<\Re x<\log r_2/r_1,\,0<\Im<2\pi\}$ modulo the translation $s\rightarrow s+2\pi i$ identification. Under the equivalence $\zeta\partial/\partial \zeta$ corresponds to $\partial/\partial s$ and $d\zeta/\zeta$ corresponds to $ds$.  
Annuli in Riemann surfaces are investigated by considering embeddings of geometric annuli.   For an inclusion of annuli $\mathcal A\subset\mathcal A'$, the modules satisfy $\mathcal M(\mathcal A)\le\mathcal M(\mathcal A')$.  A large module annulus is conformally equivalent to a geometric annulus with large ratio of boundary radii.  A sequence of annuli is degenerating provided modules tend to infinity.

For a nodal family, the total space of the family is given directly.  The standard family is the germ at the origin of the family of hyperbolas in $\mathbb C^2$.  Euclidean space $\mathbb C^2$ is a singular fibration over $\mathbb C$, given by the projection map $\pi(x,y)=x^2-y^2=t$.  The differential of the projection is $d\pi=2xdx-2ydy$.  The vector field $v=y\partial/\partial x+x\partial/\partial y$ on $\mathbb C^2$ is tangent to the fibers of $\pi$, since $d\pi(v)$ vanishes.  The vector field has constant pairing $2$ with the differential $dx/y+dy/x$.  The $0$ fiber of $\pi$ is the pair of intersecting lines $x=\pm y$.  The change of variables $z=x+y$, $w=x-y$ gives a second description of the hyperbola family.  The projection map becomes $\pi(z,w)=x^2-y^2=zw=t$ with differential $d\pi=zdw+wdz$.  The vertical vector field is $v=z\partial/\partial z-w\partial/\partial w$.  The vector field has constant pairing $2$ with the differential $dz/z-dw/w$. The $0$ fiber of $\pi$ is the union of coordinate axes.  

\section{The analytic geometry of $zw=t$.}     

For a complex manifold, we write $\mathcal O$ for the sheaf of germs of holomorphic functions, $T$ for the holomorphic tangent sheaf, $\Omega$ for the holomorphic cotangent sheaf and $K$ for the canonical bundle, the determinant line bundle of $\Omega$.  We consider the singular fibration of a neighborhood $V$ of the origin in $\mathbb C^2$ over a neighborhood $D$ of the origin in $\mathbb C$ given for $c,c'$ positive, by $V=\{|z|<c,|w|<c'\}$ with projection map $\pi(z,w)=zw=t$ to $D=\{|t|<cc'\}$. The differential of the projection $d\pi=zdw\,+\,wdz$ vanishes only at the origin; the $t=0$ fiber crosses itself at the origin.  Solving for $z$, for $t\ne 0$, the fiber of $\pi$ is $(z,w)$ with $|t|/c'<|z|<c$ and for $t=0$, the fiber is the union of discs $(z,0)$ with $|z|<c$ and $(0,w)$ with $|w|<c'$ in $\mathbb C^2$.  The family $V$ over $D$ is a family of annuli degenerating to a one point union of a $z$ disc and a $w$ disc in $\mathbb C^2$.  Alternatively, $V$ over $D$ is a germ at the origin of the family of hyperbolas limiting to the union of coordinate axes. 

\begin{figure}[htbp] 
  \centering
  \includegraphics[bb=0 0 556 556,width=2.5in,height=2.5in,keepaspectratio]{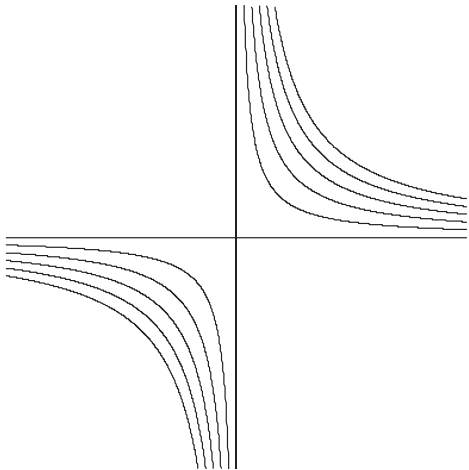}
  \caption{The family of hyperbolas $zw=t$.}
  \label{fig:hyperb}
\end{figure}

The vector field 
\[
v\,=\,z\frac{\partial}{\partial z}-w\frac{\partial}{\partial w}\in T_V,
\]
is vertical on $V-\{0\}$, since $d\pi(v)$ vanishes. Let $v'$ be another vertical vector field, non vanishing on $V-\{0\}$.  Since on $V-\{0\}$, $\ker d\pi$ is rank one, it follows that $v'=fv$, for $f$ a section of $O(V-\{0\})$.  By Hartog's Theorem \cite{Narbk}, $f$ is analytic on $V$ and since $v,v'$ are non vanishing on $V-\{0\}$, it follows that $f$ is non vanishing on $V$.  The observations provide that the condition {\em vertical vector field} defines a line bundle $\mathcal L$ over $V$ and $v$ represents a non vanishing section (the vector field $v$, a section of $\Omega_V$, vanishes at the origin; the corresponding line bundle section is non vanishing; $\mathcal L$ is not a sub bundle of $\Omega_V$).  In general, vertical vector fields that vanish at most in codimension $2$ correspond to non vanishing vertical line bundle sections. 

The $0$ fiber of $V$ over $D$ is an example of a noded Riemann surface \cite[\S 1]{Bersdeg}, alternatively a nodal curve.  The $0$ fiber is normalized by removing the origin of $\mathbb C^2$ to obtain a $z$ disc, punctured at the origin and a $w$ disc, punctured at the origin. The origins are filled in to obtain disjoint discs.  Analytic quantities on the $0$ fiber, lift to analytic quantities on the normalization.  The sheaf of regular $1$-differentials, \cite[\S 1]{Bersdeg}, equivalently the dualizing sheaf \cite[Dualizing sheaves, pg. 82]{HMbook} associates to the $0$ fiber: Abelian differentials $\beta_z,\beta_w$ with at most simple poles respectively at the origin for the $z,w$ discs and the important residue matching condition $\res \beta_z\,+\,\res\beta_w=0$.  Harris and Morrison explain for a nodal curve, the dualizing sheaf plays the role in Kodaira-Serre duality of the canonical bundle for a compact Riemann surface.

The fiber tangent spaces of $V-\{0\}$ over $D$ are subspaces of $T_V$.  The meromorphic differential
\[
\alpha\,=\,\frac{dz}{z}-\frac{dw}{w},
\]
is a functional on the fiber tangent spaces of $V-\{0\}$. For the $0$ fiber, it is immediate that $\alpha$ is a section of the dualizing sheaf.  The differential $\alpha$ satisfies the relations
\[
\alpha\wedge d\pi\,=\,2dz\wedge dw\qquad\mbox{and}\qquad\alpha(v)\,=\,2.
\]
The differential is uniquely determined modulo the submodule $\mathcal O(d\pi)\subset\Omega_V $  by each relation.  To  motivate the definition of the {\em relative dualizing sheaf} for $V$ over $D$, we consider the coset of $\alpha$ in $\Omega_V/\mathcal O(d\pi)$.  The finite pairing $\alpha(v)$ provides that the coset has holomorphic sections on $V-\{0\}$.  Let $\alpha'$ be another meromorphic differential on $V$ with $\alpha'\wedge d\pi$ holomorphic and non vanishing on $V$.  On $V-\{0\}$, the quotient sheaf $\Omega_V/\mathcal O(d\pi)$ is pointwise rank one and thus $\alpha'=f\alpha$ for $f$ holomorphic on the domain.  Again by Hartog's theorem \cite{Narbk}, the function $f$ is analytic on $V$.  The relation $\alpha'\wedge d\pi=f\alpha\wedge d\pi$ and non vanishing of the first quantity imply that $f$ is non vanishing on $V$.  In particular, the differentials $\alpha'$ with  $\alpha'\wedge d\pi$ non vanishing, considered as elements of $\Omega_V/\mathcal O(d\pi)$, define a line bundle over $V$; the line bundle is not a sub bundle of $\Omega_V$.  The differentials $\alpha'$ correspond to non vanishing holomorphic sections of the line bundle.  The constant relation $\alpha(v)=2$, shows that the line bundle is the dual of the vertical line bundle $\mathcal L$.  The relative dualizing sheaf $\omega_{V/D}$ is defined to have sections of the relative cotangent bundle $\coker(d\pi:\pi^*\Omega_D\rightarrow\Omega_V)$, given by differentials satisfying the polar divisor and residue conditions \cite[Dualizing sheaves, pg. 84]{HMbook}.  In particular, we have the sheaf equality $d\pi:\pi^*\Omega_D=\mathcal O(d\pi)$ and the quotient $\Omega_V/\mathcal O(d\pi)$ is the intended cokernel.  We have described the  relative dualizing sheaf and shown that it is dual to the vertical line bundle $\mathcal L$. 

Since the total space $V$ is smooth, there is a description of $\omega_{V/D}$ in terms of the canonical bundle $K_V$ and pullback $\pi^*K_D$ of the canonical bundle of $D$ \cite[Dualizing sheaves, pg. 84]{HMbook}; in particular
\begin{equation}\label{omegaK}
\omega_{V/D}\,\simeq\, K_V\otimes\pi^*K_D^{\vee},
\end{equation}
where $^{\vee}$ denotes the dual.  As above, a differential $\alpha'$ satisfying the polar divisor and residue conditions, determines a coset in $\Omega_V/\mathcal O(d\pi)$.  A differential $\alpha'$ and a non vanishing section $\beta$ of $\mathcal O(d\pi)$ together determine the element $\alpha'\wedge\beta\otimes\beta^{\vee}$ in $K_V\otimes\pi^*K_D^{\vee}$.  The association $\alpha'\mod\mathcal O(d\pi)\longleftrightarrow \alpha'\wedge\beta\otimes\beta^{\vee}$ is independent of the particular choice of $\beta$ and realizes the sheaf isomorphism (\ref{omegaK}).

The power $\omega_{V/D}^{k}$ of the relative dualizing sheaf, alternatively the sheaf of regular $k$-differentials \cite[\S 1]{Bersdeg}, is defined to have sections $\eta$ of 
$\Omega_V^{\otimes k}/\Omega_V^{\otimes k-1}\otimes\mathcal{O}(d\pi)$, with at most order $k$ poles at the origins of the normalization of the $0$ fiber and residues matching $\res \eta_z=(-1)^k\res\eta_w$ for the forms $\eta_z,\eta_w$ on the normalization.   The cosets are represented by $k$-fold products of sections of $\omega_{V/D}$.  The power $\omega_{V/D}^{k}$ is dual to the power $\mathcal L^k$ of the vertical line bundle.  In the next section, we give a direct description realizing the isomorphism (\ref{omegaK}). 
 
An annulus $A_{t,c,c'}=\{|t|/c'<|\zeta|<c\}$ in $\mathbb C$ maps into the $t$ fiber of $V$ by $\zeta\rightarrow(\zeta,t/\zeta)$.  The vector field $\zeta\,\partial/\partial \zeta$ on the annulus pushes forward to the vector field $v=z\,\partial/\partial z-w\,\partial/\partial w$ and the differential $\alpha=dz/z-dw/w$ pulls back to the differential $2d\zeta/\zeta$.  We introduce a boundedness notion for holomorphic differentials on annuli.  
   
\begin{defn} For positive constants $M$, $\rho_1<\rho_2<1$, a $k$-differential $\eta$ on the annulus $A_{t,c,c'}$ is {\em band bounded}, provided $\big|\eta\,(d\zeta/\zeta)^{-k}\big|\le M$ for $\zeta$ satisfying $|t|/(c'\rho_2)\le|\zeta|\le |t|/(c'\rho_1)$ and satisfying $\rho_1c\le|\zeta|\le\rho_2c$.  A sequence of $k$-differentials $\eta_t$ on annuli $A_{t,c,c'}$ with $t$ tending to zero, is band bounded provided the differentials $\eta_t$ on $A_{t,c,c'}$ are band bounded for positive constants $M$, $\rho_1$ and $\rho_2$ and all small $t$.
\end{defn}

The product $\eta\,(d\zeta/\zeta)^{-k}$ is a function with values not depending on choice of domain coordinate.  By the maximum principle, $\big|\eta\,(d\zeta/\zeta)^{-k}\big|$ is bounded on the annulus by the bounds for the bands.  For the annulus inner band and change of variable $\zeta=t/w$, the magnitude condition poses that the function of $w$ is bounded by the constant $M$.  There is also a formulation of band bounded in terms of collars for hyperbolic metrics.  In general, an annulus is described either by an embedding into a Riemann surface or an embedding into the cyclic cover corresponding to a core curve.  By Uniformization, an annulus with (complete or incomplete) hyperbolic metric is realized as a subdomain in a fiber of the standard annular family of hyperbolic metrics.  The standard annular family is described by setting $c,c'$ equal to unity in the definition of the nodal family $V$.   Each fiber of $V-\{0\}$ over $D$ has a complete hyperbolic metric \cite[Chap. 2, \S 7]{Wlcbms}, 
\begin{equation*}
\begin{split}
\label{hypexp}
dh^2_t=&\, \Bigl(\frac{|d\zeta|}{|\zeta|\log|\zeta|}\Bigr)^2\,\Bigl(\Theta\csc\Theta\Bigr)^2\mbox{\quad for\quad}\Theta=\frac{\pi\log |z|}{\log |t|},\ \frac{\pi\log |w|}{\log |t|},\\
=&\, dh_0^2\,\bigl(1+\frac13\Theta^2+\frac{1}{15}\Theta^4+\dots\bigr).
\end{split}
\end{equation*}
Collars for hyperbolic metrics are described as $|t|/\rho\le |z|\le\rho$, for suitable choices of $\rho<1$, \cite[Chap. 4]{Busbook}.  On the collar family, $(d\zeta/\zeta)(dh_t^2)^{-1}$ and its reciprocal are uniformly bounded.  The bounded magnitude condition can be posed as 
$\eta(dh_t^2)^{-k}$ has uniformly bounded magnitude on collars. 

For an annulus $A_{t,c,c'}$ in $\mathbb C$ with coordinate $\zeta$, there is a natural decomposition of holomorphic functions and for differentials $\eta$ into a sum of $\eta_++\eta_0+\eta_-$ with $\eta_+$ holomorphic inside the outer boundary, $\eta_0$ a multiple of $(d\zeta/\zeta)^k$ and $\eta_-$ holomorphic outside the inner boundary.  
The decomposition is given by the Cauchy integral for $f(\zeta)$ where $\eta=f(\zeta)(d\zeta/\zeta)^k$.  In particular for $\rho$ close to $1$, the Cauchy integral is
\begin{equation}\label{Cauchy}
f(\zeta)\,=\,\frac{1}{2\pi i}\int_{|s|=\rho c}\frac{f(s)\,ds}{s-\zeta}\ -\ \frac{1}{2\pi i}\int_{|s|=c'/\rho }\frac{f(s)\,ds}{s-\zeta}.
\end{equation} 
The first integral is holomorphic for $|\zeta|<\rho c$ and has value 
\[
f_0\,=\,\frac{1}{2\pi i}\int_{|s|=\rho c}\frac{f(s)\,ds}{s},
\] 
at the origin.  The second integral is holomorphic for $c'/\rho<|\zeta|$ and vanishes at infinity.  An equivalent description is the decomposition of the Laurant series of $f$ into the positive, zero and negative powers of $\zeta$.  The decomposition for $f$ gives the decomposition $\eta_++f_0(d\zeta/\zeta)^k+\eta_-$ for the differential $\eta$.  

\begin{lemma}\textup{Normal families.}\label{normal}  A band bounded sequence of holomorphic $k$-differentials, $\eta_t$ on the $t$ fiber of $V$, with $t$ tending to zero, has a subsequence converging uniformly on compacta on $0<|z|<\rho_2c,\,0<|w|<\rho_2c'$ to a holomorphic section of the $k^{th}$ power of the dualizing sheaf of the $0$ fiber.
\end{lemma}
\begin{proof}
For suitable $\rho$, the integration circles of (\ref{Cauchy}) $|s|=\rho c$ and $|s|=c'/\rho$ are contained in the bands and the integrand numerators are suitably bounded.  The elementary estimate for the Cauchy integral provides that the integrals are bounded as functions of $\zeta$.  The first integral is bounded holomorphic for $|\zeta|<\rho c$ with value $f_0$ at the origin.  The second integral is bounded holomorphic for $1/\rho c'<|\zeta|$ with value zero at infinity.  The conclusion follows by a standard normal families argument.
\end{proof}
We note that in the proof, the Schwarz Lemma can be applied to bound the Cauchy integrals (subtracting the value at the origin from the first) by functions vanishing at $0$ and $\infty$. 
 
\section{Sections of powers of the relative dualizing sheaf; families of regular $k$-differentials.}  

Let $\caC\stackrel{\pi}{\longrightarrow}\caB$ be a holomorphic family of noded Riemann surfaces (possibly open) over a base $\caB$, with $\caC$ and $\caB$ smooth.  We require that $\pi$ is a submersion on the complement of a codimension $2$ subset.  Smooth families of Riemann surfaces and the standard nodal family are included in the considerations.  

We present an isomorphism of sheaves.  Restricting domains as necessary, let $\beta$ be a non vanishing section of the relative dualizing sheaf $\omCB$ and $\tau$ a non vanishing section of $K_{\caB}$.  Consider the association between sections $\psi$ of $\omCB^k$ (sections of a power of the relative dualizing sheaf) and $k$-canonical forms for $\caC$, sections $\Psi$ of $K_{\caC}^k$ - the association is given by the formula
\begin{equation}\label{psiPsi}
\psi\,=\,\frac{\Psi}{(\beta\wedge\pi^*\tau)^k}\,\beta^k.
\end{equation}  
Observations are in order.  On the submersion set for $\pi$, a non vanishing section of 
$\omCB=\coker(d\pi:\pi^*\Omega_{\mathcal B}\rightarrow\Omega_{\mathcal C})$ and the pullback of a local frame for $\Omega_{\mathcal B}$ together form a frame for $\Omega_{\mathcal C}$.  
It follows that the product $\beta\wedge\pi^*\tau$ is a non vanishing section of the canonical bundle $K_{\caC}$ on the submersion set, and with the codimension $2$ condition, the product is non vanishing in general - consequently the ratio $\Psi/(\beta\wedge\pi^*\tau)^k$ is a function.  The relation can be inverted to give a formula for $\Psi$ in terms of $\psi$.  In particular, the association provides a local isomorphism of sheaves. The right hand side depends on the choice of $\tau$, but is homogeneous of degree zero in $\beta$, and so is independent of the particular choice of $\beta$. The association establishes a (twisted by $\pi^*K_{\caB}^k$) isomorphism between $\omCB^k$ and $K_{\caC}^k$.  The isomorphism 
\begin{equation}\label{omegaKC}
\omCB^k\simeq (K_{\caC}\otimes \pi^*K_{\caB}^{\vee})^k 
\end{equation}
is a general form of the isomorphism (\ref{omegaK}).   Important for our considerations, $k$-canonical forms present a local model for sections of powers of the relative dualizing sheaf, alternatively a model for families of regular $k$-differentials.

We now assume that the family $\caC\stackrel{\pi}{\longrightarrow}\caB$ can be expressed as a Cartesian product of the standard nodal family $V\stackrel{\pi}{\longrightarrow}D$ and a parameter space $S$.

\begin{lemma}\textup{Families holomorphic extension.}\label{extension}  For $V'=V-\{\pi^{-1}(0)\}$ and $D'=D-\{0\}$, let $\psi$ be a band bounded section of $\omega_{V'\times S/D'\times S}^k$ over $V'\times S$.  Then $\psi$ has a unique holomorphic extension to a section of $\omega_{V\times S/D\times S}^k$ over $\{|z|<\rho_2c,|w|<\rho_2c'\}\times S$. 
\end{lemma} 
\begin{proof} For the Cauchy integral (\ref{Cauchy}), write $\psi=f\alpha^k$ for $f=f(z,w,s)$.   The Cauchy integral considerations show that the function $f$ is bounded on $V'\cap\{|z|<\rho_2c,|w|<\rho_2c'\}\times S$.  By formula (\ref{psiPsi}), using $\alpha=dz/z-dw/w$ as the reference section, the associated section of $K_{V\times S}^k$ is $\Psi=\psi\,\alpha^{-k}(\alpha\wedge \pi^*dt\wedge ds)^k$, $ds$ an $S$ canonical form.  Since $\psi\,\alpha^{-k}=f$ and $\alpha\wedge\pi^*dt=2dz\wedge dw$, we have that $\Psi=f\,(2dz\wedge dw\wedge ds)^k$.  By the Riemann extension theorem \cite{Narbk}, $f$ has a unique holomorphic extension to $V\times S$.  The formula $\psi=f \alpha^k$ gives the desired extension of $\psi$. 
\end{proof} 

Observations are in order.  The lemma provides that a $k$-differential is given on $\{|z|<\rho_2c,|w|<\rho_2c'\}\times S$ as
\[
\psi\,=\,\mathbf f(z,w,s)\,(\frac{dz}{z}-\frac{dw}{w})^k,
\]
for $\mathbf f(z,w,s)$ a holomorphic function.  The differential is given as
\[
\psi\,=\,\mathbf f(\zeta,t/\zeta,s)\,(2\frac{d\zeta}{\zeta})^k,
\]
on annuli $\{|t|/c'<|\zeta|<c\}$ mapped into fibers.   On the locus $w=0$, $dw$ vanishes and $\psi$ is given as $\mathbf f(z,0,s)(dz/z)^k$ and on the locus $z=0$, $dz$ vanishes and $\psi$ is given as $\mathbf f(0,w,s)(-dw/w)^k$.  Lemmas \ref{normal} and \ref{extension} show the role of the relative dualizing sheaf in considering limits of band bounded holomorphic $k$-differentials - the sheaf naturally appears.  More generally, if $\psi$ is bounded on annuli outer bands by $M$ and on annuli inner bands by $M|t|^{-m}$, then $\psi\zeta^m$ is band bounded.  It follows that
\[
\psi\,=\,\mathbf f(z,w,s)z^{-m}\,(\frac{dz}{z}-\frac{dw}{w})^k,
\]
for $\mathbf f(z,w,s)$ a holomorphic function.  On the locus $w=0$, $\psi$ is given as $\mathbf f(z,0,s)z^{-m}(dz/z)^k$ and the locus $z=0$ is a polar divisor provided $\psi$ is unbounded on  annuli inner bands.

The isomorphism (\ref{psiPsi}) provides a setting for understanding the zero and polar divisors of sections of the relative dualizing sheaf, alternatively holomorphic families of regular $k$-differentials.  The reference section $\beta$ is non vanishing; $\psi$ and $\Psi$ have the same divisors.  We consider two examples, using $\alpha$ as reference section and omitting the general parameter space.  First consider an example of the $w$-axis contained in the divisor 
of $\Psi$, in particular $\Psi=z\,\mathbf f(z,w) (2dz\wedge dw)^k$.  The differential $\psi$ is given as $\psi=z\,\mathbf f(z,w)(dz/z-dw/w)^k$; $\psi$ is trivial on the $w$-axis and $\psi=\mathbf f(z,0)z^{1-k}(dz)^k$ on the $z$-axis.  A second example is for the multiplicity $m$ divisor $(az+bw)^m=0,\,a\ne 0,b\ne 0,$ containing the node, in particular $\Psi=(az+bw)^m\,\mathbf f(z,w) (2dz\wedge dw)^k$ with $\mathbf f$ not vanishing at the origin.  The intersection of the divisor and a $t$ fiber, $t\ne 0$, are the solutions of $((az^2+bt)/z)^m=0$.  The differential $\psi$ evaluates to $b^mw^m\,\mathbf f(0,w)(-dw/w)^k$ on the $w$-axis and to $a^mz^m\,\mathbf f(z,0)(dz/z)^k$ on the $z$-axis.    The differential $\psi$ has $2m$ zeros on each collar and an order $m$ zero as a section of $\mathcal O((K(o))^k)$ on each branch of the normalization of the $0$ fiber (the origin determines the point divisor $o$ on each branch of the normalization).  In general, if $\mathbf f(z,w)$ is non vanishing on bands, then the winding principle can be applied to $\mathbf f(z,t/z),\,\mathbf f(z,0)$ and $\mathbf f(0,w)$ - to find that with the above approach for the nodal fiber, the count of zeros on fibers is constant.  

There is a direct relationship of the isomorphism (\ref{omegaKC}), the isomorphism (\ref{psiPsi}) and the Poincar\'{e} residue map for adjunction (see \cite[Adjunction Formula II, pg. 147]{GHbook}).  The Poincar\'{e} residue $K_{\caC}(F)\stackrel{P.\,R.}{\longrightarrow}\omCB\big|_F$ for a fiber $F$ of $\caC$ over $\caB$ is given as follows. In general, for a divisor $\caD\subset M$ in a complex manifold $M$ given as a locus $h(z)=0$, then the Poincar\'{e} residue $K_M(\caD)\longrightarrow K_{\caD}$ is given by the association of $\Phi$, a section of $K_M(\caD)$ with polar divisor $\caD$, to $\phi$ a section of $K_{\caD}$, provided $\Phi=dh/h\,\wedge\phi$.  

Consider the iterated Poincar\'{e} residue for a family $\caC\stackrel{\pi}{\longrightarrow}\caB$, with $\caC$ and $\caB$ smooth and $\operatorname{det}d\pi$ vanishing at most in codimension $2$.  Let $z=(z_1,\dots,z_{n+1})$ be a local coordinate for $\caC$ and $t=(t_1,\dots,t_n)$ a local coordinate for $\caB$.  Let the projection be given as 
$t=\pi(z)=(\pi_1(z),\dots,\pi_n(z))$.   Given a particular value for $t$, the successive  divisions by $(\pi_j(z)-t_j)$ and computing Poincar\'{e} residues: for the divisor $\pi_n(z)=t_n$ inside $\caC$; for the divisor $\pi_{n-1}(z)=t_{n-1}$ inside $\{\pi_n=t_n\}$; $\ \dots\ $; for the divisor $\pi_1(z)=t_1$ inside $\{\pi_n=t_n,\pi_{n-1}=t_{n-1},\dots,\pi_2=t_2\}$ begins 
with $\Phi$ a section of $K_{\caC}$, and assigns a $1$-form $\phi$ a section of $\Omega_{\caC}$, satisfying
\begin{equation}\label{residue}
\frac{\Phi}{\prod_{j=1}^n(\pi_j(z)-t_j)}\,=\,\frac{d\pi_1}{(\pi_1(z)-t_1)}\wedge\cdots\wedge \frac{d\pi_n}{(\pi_n(z)-t_n)}\wedge\phi,
\end{equation}
or clearing denominators, gives the equation $\Phi=d\pi_1\wedge\cdots\wedge d\pi_n\wedge\phi$. Since $d\pi_1\wedge\cdots\wedge d\pi_n=(d\pi)^*(dt_1\wedge\cdots\wedge dt_n)$, $\phi$ is determined modulo $\pi^*K_{\caB}$ and equation (\ref{residue}) is a form of the isomorphism (\ref{omegaKC}). Equivalently, the iterated Poincar\'{e} residue map for fibers of $\pi$ is equivalent to the isomorphism (\ref{psiPsi}).  In Masur's paper \cite{Msext} on extension of the Weil-Petersson metric, the iterated Poincar\'{e} residue map is used as the context in Proposition 4.2 for the construction of holomorphic families of $1$-differentials from canonical forms. Formula (\ref{psiPsi}) provides a direct description of the association with canonical forms. 

A basic matter is the construction of holomorphic families of $k$-differentials for families of possibly noded Riemann surfaces; the matter includes showing that $k$-differentials form holomorphic vector bundles over deformation spaces.  In Bers' original article \cite{Bersdeg}, `Spaces of degenerating Riemann surfaces', the notion of regular $k$-differentials (sections of the dualizing sheaves of individual curves) is formulated, but a notion of a {\em holomorphic family} of regular $k$-differentials (sections of the relative dualizing sheaf) is not formulated.  Bers parameterizes families of possibly noded Riemann surfaces in terms of families of Kleinian groups, groups that vary holomorphically in $SL(2;\mathbb C)$.  He uses Poincar\'{e} series of rational functions to describe local families of regular $k$-differentials that vary holomorphically on compacta within the domains of discontinuity.  He explains that the families generically give bases for regular $k$-differentials \cite[\S 4]{Bersdeg}.  A next step is to show that the bases give local frames for vector bundles of differentials.  For a suitable family of Kleinian groups, Poincar\'{e} series constructions can give families of differentials varying holomorphically and converging uniformly on compacta  in domains of discontinuity.  Uniform convergence on compacta provides the band bounded condition for annuli in the quotient.  In this setting, Lemma \ref{extension} can provide that families of Poincar\'{e} series give sections of powers of relative dualizing sheaves. Formula (\ref{psiPsi}) provides local expansions of sections.


\providecommand\WlName[1]{#1}\providecommand\WpName[1]{#1}\providecommand\Wl{W%
lf}\providecommand\Wp{Wlp}\def\cprime{$'$}

\end{document}